\documentclass[12pt]{article}
\usepackage{row-article}
\usepackage{graphicx}
\usepackage[utf8]{inputenc}
\usepackage{eufrak}
\usepackage{makeidx}
\begin{document}
\title{Compact Complex Surfaces with No Nonconstant Meromorphic Functions %
}
\author{Raymond O. Wells, Jr.
\footnote{Jacobs University Bremen; University of Colorado at Boulder; rwells@colordo.edu} }       
\maketitle

\begin{abstract}
In 1949 Siegel gave an example of a complex two-torus with no nonconstant meromorphic functions. In 1964 Kodaira showed that compact complex surfaces with no nonconstant meromorphic functions must be of the following three types: tori, Hopf type surfaces with first Betti number equal to one, and K3 surfaces. In this paper we show that surfaces of these three types have a dense set of surfaces in their natural moduli spaces with no nonconstant meromorphic functions.
\end{abstract}

\section{Introduction}
Let \(X\) be a connected compact complex manifold, and let \(\SM(X)\) be the field of meromorphic functions on \(X\). We let \(\deg X\) denote the {\em transcendence degree} of \(\SM(X)\) over the field \(\BC\) of complex numbers, which we here identify with the field of constant complex-valued functions on \(X\).

In 1955 Siegel  \cite{siegel1955} showed that 
\be
\label{eqn:siegel1955}
0\le \deg X \le n, \;\textrm{where}\;n=\dim X.
\ee
This result has a long history going back to the first efforts by Riemann and Weierstrass to show that there are not \(n+1\) algebraically independent Abelian function on \(\BC^n\).  In his paper Siegel has a very informative history of the results and proofs leading up to (\ref{eqn:siegel1955}).

We note that when Riemann first introduced Riemann surfaces in 1857 \cite{riemann1857}, one of the primary questions he addressed was the question of the existence of meromorphic functions. In this initial paper he did show the existence of nonconstant meromorphic functions on a Riemann surface, this then leading up to the Riemann-Roch theorem. This was the context of a Riemann surface being a branched covering of the compactified complex plane. When Weyl  introduced for the first time an abstract Riemann surface \cite{weyl1913}, he also showed, following the original outline of Riemann, that any compact Riemann surface has nonconstant meromorphic functions, and that the Riemann-Roch theorem was valid as well. Thus, any connected compact complex manifold \(X\) of dimension one has nonconstant meromorphic functions. We shall see that for higher dimensions, this is no longer the case.

In this paper we want to concentrate on connected compact complex surfaces, i.e., complex manifolds of dimension two, and we will call them simply {\em surfaces} (assumed compact and connected).  Thus, if \(X\) is a surface, we have
\[
0\le \deg(X)\le 2.
\]
In one of Kodaira's seminal papers  ``On the structure of compact complex analytic surfaces. I'' \cite{kodaira1964}, he shows that:
 \begin{enumerate}
 \item If \(\deg X =2\), then \(X\) must be a projective algebraic manifold%
 \footnote{In the general case, if \(\deg X = n\), where \(X\) is \(n\)-dimensional, then \(X\) is not necessarily a projective algebraic manifold; instead it is referred to as a Moishezon manifold, which is always bimerorphically equivalent to a projective algebraic manifold.  See Shafarevich \cite{shafarevich1997} for a discussion of this topic.}%
 .
 
 \item If \(\deg X = 1\), then \(X\) is a type of {\em elliptic surface}, a fiber space of elliptic curves (with a finite number of singular curves), parametrized by an algebraic curve. 
 \end{enumerate}
 
 What can one say about surfaces of degree zero? Siegel gave an example in his 1949 Lecture Notes \cite{siegel1949} of a complex two-torus of degree zero, which we will see again later in this paper.  In 1964 Kodaira \cite{kodaira1964} showed that there are topological restrictions on their structure. We let \(b_i=b_i(X)\) be the Betti numbers of a surface \(X\) and let \(K_X=\wedge^2 T^*(X)\) be the canonical bundle of \(X\), then Kodaira showed that if a surface \(X\) has no exceptional curves%
\footnote{An \emph{exceptional curve} in a surface \(X\) is a curve with self-intersection \(c\cdot c =-1\); Grauert \cite{grauert1962} showed that any surface \(X\) with exceptional curves \(c_1,\ldots,c_k\) is biholomorphic to a surface \(\tilde X\) with with a finite number of quadratic transforms (blowups) at points \(p_1,\ldots,\p_k\), where the curves \(c_i\) are the quadratic transforms of the points \(p_i\) in \(\tilde X\).}
then there are three possibilities:
 \begin{enumerate}
 \item \(b_1=4\), and \(X\) is a complex torus.
 
 \item \(b_1=1\), and \(X\) is a  Hopf surface or other similar non-Kähler surfaces (Inoue, Enriques-Hirzebruch, etc.).
 
  \item \(b_1 = 0\), the canonical bundle \(K_X\) of \(X\) is trivial, and \(X\) is a K3 surface.
 
 \end{enumerate}
 
In Shafarevich's book \cite{shafarevich1997} there is a simple criterion for a  Hopf surface to be of degree zero which immediately yield examples of degree zero Hopf surfaces.  We did not find in the literature any examples of K3 surfaces of degree zero. In particular, in Kodaira's paper \cite{kodaira1964} he only shows that degree zero surfaces must be of the three types above, not that such manifolds exist with this property. 
 
 However, these examples suffice to show that the theory of several complex variables again distinguishes itself from function theory of one variable, as has been the case since Hartogs proved his fundamental theorem in 1906 \cite{hartogs1906} concerning simultaneous analytic continuation across compact subsets of \(\BC^2\).
 
 In this paper we show that for each of these three cases: complex two-tori, Hopf surfaces, and K3 surfaces,  the set of degree zero surfaces is {\em dense} in the natural moduli spaces for each of these classes of surfaces. We don't consider the other non-Kähler surfaces here, but their treatment would, in principle, be similar to our study of Hopf surfaces, as they are defined in a similar manner as a quotient of a specific type of domain in \(\BC^2\) by a discrete group.
 
 To be more specific in the first case, we consider complex two-tori defined by period matrices of type
 \be
 \label{eqn:period-matrix}
 \O=(I,Z),
 \ee
 where \(I\) is the \(2\times 2\) identity matrix and \(Z\) is a \(2\times 2\) complex-valued matrix with \(\Im Z > 0\) (positive definite). Let
 \[
 M=\{Z: Z \;\textrm{is a}\;2\times 2 \;\textrm{complex-valued matrix with}\; \Im Z>0\}.
 \] 
 This is a moduli space of dimension four for the complex tori defined by the period matrices (\ref{eqn:period-matrix}) (it is a complete and effectively parametrized moduli space as shown by Kodaira and Spencer \cite{kodaira-spencer1958} (see Theorem 14.3 in \cite{kodaira-spencer1958}). 
 
 Our result in this case is that there is a dense set of points \(M_0\subset M\) such that%
 \footnote{In fact, \(M_0\) is a set of second category in the four-dimensional metric space \(M\).}%
  for each \(Z\in M_0\), the torus \(T_Z\) defined by the period matrix \(\O=(I,Z)\) has degree xero.  The other two cases are formulated similarly, and we will consider the detailed results for all three cases in the following three sections.

\section{Two-dimensional tori}
Let \(M_2\) denote the vector space of \(2\times 2\) complex-valued matrices, which we will denote generically by
\[
Z=\left(\ba{cc}
	z_{11}&z_{12}\\
	z_{21}&z_{22}\ea\right).
\]
We let
\[
M:=\{Z\in M_2: \Im Z>0\},
\]
where \(\Im Z\) denotes the imaginary parts of the coefficients of \(Z\) and where \(\Im Z>0\) means that  the matrix \(\Im M\) is positive definite. 

Let \(I\) be the identity matrix in \(M\), and let
\be
\label{eqn:period-matrix}
\O = (I,Z), Z\in M,
\ee
which we will call a {\em normalized period matrix}. Let \(\{\om_1,\om_2,\om_3,\om_4\}\) be the columns of \(\O\), and we see that these four vectors in \(\BC^2\) are linearly independent over the real numbers.  Let \(\L\) be the lattice generated by these vectors, i.e., linear combinations of of the form
\[
m_1\om_1+ m_2\om_2 + m_3\om_3 +m_4\om_4 \;\;m_i\in \BZ,
\]
and let 
\[
T_Z := \BC^2/\L,
\]
which is a complex torus of two dimensions defined by the period matrix \(\O\).

Let \(Z\in M\), then an {\em Abelian function} \(f\) on \(\BC^2\) with respect to the period matrix
\[
\O=(I,Z)=(\om_1,\om_2,\om_3,\om_4),
\]
is a meromorphic function on \(\BC^2\) which satisfies
\[
f(z+\om_j)=f(z), \;j=1,\ldots,4.
\]

A {\em degenerate Abelian function} is an Abelian function \(f(z)=f(z_1,z_2)\) which, after a linear change of variables in \(\BC^2\) 
\[
\left(\ba{c}
	\z_1\\
	\z_2\ea\right) =\left(\ba{cc}
		\a&\b\\
		\g&\d\ea\right)\left(\ba{c}
			z_1\\
			z_2\ea\right), \; \a\d-\g\b\ne0,
\]
becomes a a nonconstant elliptic function of one complex variable, i.e.,
\[
\tilde f(\z_1,\z_2)= f(z_1(\z_1,\z_2),z_2(\z_1,\z_2))=\tilde f(\z_1),
\]
where \(\tilde f(\z_1)\) is an elliptic function. 

If \(\tilde f\) is not a constant, then it is an elliptic function in the \(\z_1\) complex plane which is  also periodic with respect to the four complex numbers which are the first row of the transformed period matrix:
\[
\left(\ba{cc}
	\a&\b\\
	\g&\d\ea\right)\O,
\]
which have the form:
\[
(\a,\b,\a z_{11}+\b z_{21}, \a z_{12}+ \b z_{22}).
\]
Since \(\tilde f(\z_1)\) is an elliptic function, it has two independent periods \(\t_1\) and \(\t_2\) which generate all the periods of \(\tilde f\) in the complex \(\z_1\)-plane.  This means that there are integers \(p_1,p_2,q_1,q_2,r_1,r_2,s_1,s_2\) which satisfy
\be
\label{eqn:one-variable-periods}
\ba{rcl}\a&=&p_1\t_1 +p_2 \t_2,\\
\b&=&q_1\t_1+q_2 \t_2,\\
\a z_{11} +\b z_{12}&=&r_1 \t_1 +r_2 \t_2,\\
\a z_{21} +\b z_{22}&=&s_1 \t_1 +s_2 \t_2.
\ea
\ee
We will return to this in our analysis below.

An Abelian function \(f\) on \(\BC^2\) is {\em nondegenerate} if it is not a degenerate Abelian function. See Siegel \cite{siegel1949} or Conforto \cite{conforto1956} for an overview of nondegenerate and degenerate Abelian functions on \(\BC^n\).  We will use some of the results from these books in our discussion below.

Let 
\be
\label{eqn:general-period-matrix}
\O= (\om_1,\ldots,\om_{2n}),
\ee
be a general period matrix in \(\BC^2\), not necessarily the more specialized normalized period matrices we considered above in (\ref{eqn:period-matrix}). Here \(\{\om_1,\ldots,\om_{2n}\}\) are \(2n\) vectors in \(\BC^n\) which are linearly-independent over the real numbers, and we let \(T_\O\) denote the complex torus corresponding to \(\O\). A period matrix \(\O\) of the form (\ref{eqn:general-period-matrix}) is said to be a {\em Riemann matrix} if there exists a skew-symmetric \(2n\times 2n\) matrix \(P\) with coefficients which are rational numbers such that
\be
\label{eqn:riemann-matrix}
\ba{rcl}\O P \O^t &=& 0,\\
-i\O P \overline \O^t &>&0.
\ea
\ee
Here \(A^t\) denote the transpose of a matrix \(A\), and \(\overline A\) denotes the complex conjugate of the entries of a matrix \(A\).  

An important theorem concerning complex tori is that a complex torus of \(n\)-dimensions is a projective algebraic manifold if and only if the period matrix \(\O\) defining the torus is a Riemann matrix (this is a consequence of the Kodaira embedding theorem; see, e.g., Wells \cite{wells2008}). A different and related result is that for a given period matrix \(\O\), the matrix \(\O\) is a Riemann matrix if and only if there exists a nondegenerate Abelian function on \(\BC^n\) with respect to 
\(\O\). See Siegel \cite{siegel1949} or Conforto \cite{conforto1956} for a proof of this. It follows from this that 
\[
\deg T_\O =n \;\textrm{if and only if}\; \O\; \textrm{is a Riemann matrix}.
\]

In our two-dimensional case of normalized period matrices, we see now that there are three cases to consider. First, \(\deg T_Z=2\) if and only if \((I,Z)\) is a Riemann matrix.  Secondly, if \(\deg T_Z=1\), then there must be a nonconstant degenerate Abelian function \(f\) on \(T_Z\).  The third case is that \(\deg T_Z=0\), i.e., there are no nonconstant meromorphic functions on \(T_Z\). We will give explicit criteria for the first two cases in terms of the matrices \(Z\in M\), which will then determine the nature of the set of matrices \(Z\in M\) for which \(\deg T_Z=0\).

We now consider period matrices for \(n\)-dimensional complex tori, and we will restrict ourselves again to the two-dimensional case somewhat later.  Two period matrices \(\O_1\) and \(\O_2\) are said to be \emph{equivalent} if there is a nonsingular matrix \(C\in M_k(\BC)\) and a nonsingular matrix \(N\in M_{2n}(\BZ)\) such that
\[
\O_1=C\O_2N.
\]
Two equivalent period matrices generate tori \(T_{\O_1}\)  and \(T_{\O_2}\) which are biholomorphic to each other. Using this equivalence relation, a variety of canonical forms for Riemann matrices were formulated in the late nineteenth and early twentieth century.  These are summarized quite completely in Conforto's book \cite{conforto1956} (see, in particular, the table on p. 90). 
.  

Coming back to the two-dimensional case, we see from Conforto's book that the the set of all Riemann matrices are equivalent to the canonical Riemann matrices of the form
\[
\O_n := (I_n, Z), I_n=\left(\ba{cc}
	1&0\\
	0&n
	\ea\right)
	, \;\textrm{for}\;Z\in M, Z=Z^t,
\]
for all positive integers \(n\). We note that, except for \(n=1\), these are not normalized period matrices. However, by left multiplication by \(I_n^{-1}\), we obtain the equivalent family of all normalized Riemann matrices of the form
\[
\O=(I,I_n^{-1}Z),
\]
for \(Z\in M\) and where \(Z\) is symmetric.  Therefore we set
\[
S_n:= \{Z\in M: z_{21}=nz_{12},\; n= 1,2,\dots\}
\]
and 
\[
\{\O=(I,Z), Z\in S_n\}, \; n=1,2,\ldots,
\]
consists of all Riemann matrices in normalized form.

We  let
\[
S:=\bigcup_{n=1}^\infty S_n.
\]
Each \(S_n\) is a linear hyperplane in the vector space \(M_2(\BC)\) intersecting the open set \(M\subset M_2(\BC)\), and, as such, it is a closed subset of \(M\) with no interior points, thus \(S\) is a set of first category in \(M\). Therefore we obtain that \(T_Z\) has degree two if and only if \(Z\in S\). This is a useful criterion for complex tori of degree two.

We now want to give a similar criterion for complex two-tori of degree one. We suppose now that \(Z\in M\) and \(\deg T_Z=1\). Thus there is a nonconstant degenerate Abelian function \(f\) on \(X_Z\).  It follows from our earlier discussion that there is a change of variables of the form 
\[
C=\left(
\ba{cc}
\a&\b\\
\g&\d\ea
\right),
\]
and periods \(\t_1,\t_2\in \BC\) with \(\Im \t_1/\t_2\ne 0\) and integers \(p_1,p_2, q_1,q_2, r_1,r_2,s_1,s_2\) satisfying (\ref{eqn:one-variable-periods}).  Moreover, the matrix
\[
\left(
\ba{cc}
p_1&p_2\\
q_1&q_2\\
r_1&r_2\\
s_1&s_2
\ea
\right) 
\]
has maximal rank.

Eliminating \(\a\) and \(\b\) from the four equations in (\ref{eqn:one-variable-periods}) we find that
\bean
(p_1\t_1+p_2\t_2)z_{11}+(q_1\t_1+q_2\t_2)z_{21}&=&r_1\t_1+r_2\t_2\\
(p_1\t_1 +p_2\t_2)z_{12}+(q_1\t_1+q_2\t_2)z_{22}&=&s_1\t_1+s_2\t_2,
\eean
which gives
\bean
(p_1z_{11}+q_1z_{21}-r_1)\t_1 +(p_2z_{11}+q_2z_{21}-r_2)\t_2&=&0,\\
(p_1z_{12}+q_1z_{22}-s_1)\t_1+ (p_2z_{12}+q_2z_{22}-s_2)\t_2&=&0.
\eean
Since \(\t_1\) and \(\t_2\) are linearly independent over the real numbers, it follows that the determinant of this linear system must be zero. Hence
\be
\label{eqn:determinant}
\ba{c}(r_1s_2-s_1r_2)+(s_1p_2-p_1s_2)z_{11}+ (p_1r_2-r_1p_2)z_{12}+(s_1q_2-q_1s_2)z_{21}\\ 
\mbox{\hspace{.5in}} +(q_1r_2-r_1q_2)z_{22} + (p_1q_2-q_1p_2)(z_{11}z_{22} -z_{12}z_{21})=0
\ea
\ee

Let \(m=(m_0,m_1,m_2, m_3,m_4,m_5)\) be a sextuple of integers, and assume that \(m_i=\ne0\) for some \(i, 1\le i \le5\), which we call an \emph{admissible sextuple}.  Then, for such an admissible sextuple,  we define the algebraic variety \(R_m\) in \(M\) by the equation
\be
\label{eqn:variety}
\ba{c}
R_m:=\{ Z\in M: m_0 +m_1z_{11} +m_2z_{12} +m_3z_{21}+ m_4z_{22} \\
\mbox{\hspace{.5in}}+ m_5(z_{11}z_{22}-z_{12}z_{21})\}.
\ea
\ee
We see that \(R_m\) is of codimension one. 

We note that the coefficients of the polynomial in (\ref{eqn:determinant}) are precisely the minors of the \(4\times 2\) matrix
\be
\label{eqn:period-coefficients}
\left(\ba{cc}
p_1&p_2\\
q_1&q_2\\
r_1&r_2\\
s_1&s_2\ea
\right)
\ee
Now we are assume that \(Z\in M\) and that \(\deg T_Z=1\), then it follows that \(Z\) is a point on \(R_m\).  Namely, at least one of the minors must be nonzero.  If the minor \(r_1s_2-s_1r_2\) is the only nonzero minor, then equation (\ref{eqn:determinant}) couldn't be satisfied, and so at least one other minor must also be nonzero, and thus \(Z\in R_m\) for some suitable admissible \(m\). 

Let now
\[
R= \bigcup_{m}R_m,
\]
where we sum over admissible sextets \(m\). Then any \(Z\in M\) with \(X_Z\) of degree one must be a point in \(R\), and as we saw earlier, any \(Z\in M\) with \(X_Z\) of degree two must be a point in \(S\), so we have the theorem.
\bthm
Let \(X_Z\) be a  two-torus with the period matrix \((I,Z)\), then, if \(\deg X_Z\ge 1\), then \(Z\in S\cup R\).
\ethm
Since \(S\cup R\) is a countable union of subvarieties of \(M\), each of codimension one, we have the following theorem. Let 
\[
M_0:=M- (S\cup R).
\]
\bthm
\label{thm:tori}
\(M_0\) is of second category in \(M\), and in particular is dense, and for each \(Z\in M_0\),
\[
\deg T_Z=0.
\]
\ethm
 
 We will close this section with a discussion of two examples. As mentioned in the Introduction, Siegel gave in his lecture notes \cite{siegel1949}  an example of a two-torus of degree zero which we describe now. Let
 \[
 \tilde \O_1=\left(\ba{cccc}
 1&0&\sqrt{-2}&\sqrt{-5}\\
 0&1&\sqrt{-3}&\sqrt{-7}
 \ea\right),
 \]
This is Siegel's example, but we need to choose a \emph{sign} for the square roots (other choices work as well), and we can reorder the vectors slightly, so we let 
\[
 \O_1=\left(\ba{cccc}
 1&0&\sqrt{5}i&\sqrt{2}i\\
 0&1&\sqrt{7}i&\sqrt{3}i
 \ea\right),
 \]
 and we let
 \[
 Z_1=\left(\ba{cc}
 \sqrt{5}i&\sqrt{2}i\\
 \sqrt{7}i&\sqrt{3}i\ea\right),
 \]
 which we see is a point in the moduli space \(M\).  We see immediately that \(\O_1\) is \emph{not} a Riemann matrix since \(\sqrt{7}i\) is \emph{not} of the form \(n\sqrt{2}i\), for some positive integer \(n\). 
 
 Suppose that \(Z\in R\). then there is an admissible sextuple \(m\) such that
 \[
 m_0+m_1(\sqrt{2}i) + m_2(\sqrt{5}i) +m_3(\sqrt{3}i)+m_4(\sqrt{7}i)+m_5(\sqrt{15}-\sqrt{14})=0.
 \]
This means that
\bean
m_0+m_5(\sqrt{15}-\sqrt{14})&=&0,\\
m_1\sqrt{2}+m_2\sqrt{5}+m_3\sqrt{3} +m_4\sqrt{7}&=& 0.
\eean
It is easy to conclude from these two equations that \(m_i=0, i=0,\ldots,5\), which contradicts \(m\) being an admissible sextuple.  Hence \(\Z_1\) is not in either \(S\) or \(R\), and must therefore correspond to a torus \(T_1\) which is of degree zero.

Let us now look at Shafarevich's example of a two-torus of degree zero (in his book \cite{shafarevich1997}).  We let
\[
\O_2=\left(\ba{cccc}
1&0&i&\sqrt{2}\\
0&1&0&i\ea\right),
\]
and we let 
\[
Z_2=\left(\ba{cc}
i&\sqrt{2}\\
0&i\ea\right),\]
which is also a point in \(M\).

We see immediately that this is not a Riemann matrix. However, \(Z_2\) is a point of \(R\) for some admissible sextuplet \(m\).  Namely, suppose that \(m_1=-m_4=n\), where \(n\) is a nonzero integer, and we set 
\[
m_0=m_2=m_3=m_5=0,
\]
then we find that
\bean
m_0+ m_1(i)+m_2(\sqrt{2}) +m_3(0)+m_4(i) +m_5(-1)&=&ni-ni,\\
&=&0.
\eean
Thus 
\[
Z_2 \in R_{0,-n,0,0,n,0},
\]
and is not a point in \(M_0\), whereas Siegel's example was.

This shows that the dense subset \(M_0\) of Theorem \ref{thm:tori} does not contain all the complex two-tori of degree zero.  One would need a more precise characterization of two-tori of degree one for this, as we do have for two-tori of degree two. We note that the integers \(m_i\) in (\ref{eqn:variety}) should be quadratic forms in terms of the integers \(p_i,q_j,r_k,s_l\) in (\ref{eqn:period-coefficients}), which is a type of restriction on the \(m_i\), which we have not used in our analysis, but could play a role in finding a better characterization of degree one tori.

In Section \ref{sec:K3} we will study degree zero K3 surfaces, but we note here that Kodaira showed in \cite{kodaira1964} that in the moduli space for K3 surfaces that we will use, projective algebraic K3 surfaces are dense in the moduli space, similar to the density of degree zero tori that we described above. If we let
\[
S_0:=\{Z\in M:z_{12}=0\},
\]
and set 
\[
\tilde S=S_0 \cup S,
\]
then \(\tilde S\) is a \emph{closed} subset of \(M\), and its complement \(\tilde M \subset M_0\) is an open subset of \(M\).  It follows that projective algebraic tori are not dense iin \(M\), in contrast to Kodaira's result for K3 surfaces.

\section{Hopf surfaces}
In 1948 Heinz Hopf introduced \cite{hopf1948} a class of compact complex surfaces that have become known as \emph{Hopf surfaces}. These were the first examples of compact complex manifolds which were not Kähler manifolds. Any compact Kähler manifold \(X\) has the property that the odd Betti numbers of \(X\), \(b_1, b_3,\ldots,\)  are all \emph{even} integers, and the Hopf surfaces have \(b_1=1\), as we will see below, and hence are not Kähler (see, e.g.,  Wells \cite{wells2008} as a reference for the theory of Kähler manifolds).

Let us start with a simple example, as in Hopf's paper. Let
\[
W:= \BC^2-\{0\},
\]
be the \emph{punctured complex 2-plane}, which will play an important role throughout this section. This is a simply-connected noncompact complex manifold which will play the role of a universal covering space for many of our compact complex manifolds we will be discussing in the paragraphs below. Let
\[
\g:W \rightarrow W,
\]
be a holomorphic mapping defined by 
\[
\g(z_1,z_2)=(2z_1,2z_2),
\]
and let \(\g^n\) denote \(n\) iterations of this mapping. This generates a discrete transformation group \(\G\) which acts on \(W\), and we define the quotient space
\be
\label{eqn:hopf-1}
X:=W/\G.
\ee
One can easily verify (see, e.g., Wells \cite{wells2008}, p. 200) that \(X\) is a compact complex manifold which is diffeomorphic to \(S^1\times S^3\), where \(S^1\) and \(S^3\) are the one-sphere and three-sphere, respectively, defined by
\[
S^n:=\{x\in \BR^{n+1}: x_1^2+\cdots+x_{n+1}^2=1\}.
\]
Thus \(X\) has Betti number \(b_1=1\) and is not a Kähler surface, as mentioned above.

In 1958 Kodaira and Spencer introduced their fundamental theory of deformations of compact complex manifolds \cite{kodaira-spencer1958}.  They applied this theory to a variety of examples of specific classes of complex manifolds, including tori as we discussed in the previous section, hypersurfaces of complex projective spaces, Hopf surfaces, etc. We will summarize their deformation theory for Hopf surfaces to give us a suitable moduli space with which we can formulate our theorem concerning degree zero Hopf surfaces.

Let \(M\) be the four-dimensional complex manifold of \(2\times 2\) matrices
\[
t=\left(\ba{cc}
	\a&\b\\
	\g&\d\ea\right),
\]
which are nonsingular linear transformations
\[
t:\BC^2\rightarrow \BC^2,
\]
and which have the property that 
\[
|\a+\d|>3, |(\a-\d)^2+4\b\g|<1.
\]
Then the eigenvalues of \(t\in M\) are
\[
\s\pm\sqrt{\D}, \;\textrm{where}\; \s=\frac{1}{2}(\a+\d), \D=\frac{1}{4}(\a-\d)^2+\b\g.
\]

Letting \(W=\BC^2-\{0\}\), as before, we define the holomorphic automorphisms of \(W\times M\) by
\[
\eta:(z,t)\mapsto (tz,t).
\]
Since the eigenvalues of \(t\) satisfy \(|\s\pm\sqrt{\D}|>1\), it follows that 
\[
G=\{\eta^m, m\in \BZ\}
\]
is an infinite cyclic group which is a properly discontinuous group of biholomorphic mappings with no fixed points, and thus the quotient space
\[
X:= (W\times M)/G
\]
is a  complex manifold.

Let 
\[
p:W\times M\rightarrow X
\]
be the canonical projection and thus there is  commutative diagram
\[
\ba{ccccc}
\;\;W\times M&&\stackrel{p}{\longrightarrow} &&X\\
&\searrow&&\stackrel{\p}{\swarrow}&\\
&&M&&\ea
\]
and the triple \((X, \p,M)\) is a complex-analytic family of complex manifolds. We let
\[
X_t:=\p^{-1}(t)= W/G_t,
\]
where \(G_t\) is the infinite cyclic group generated by \(t\) acting on \(W\). The manifold \(X_t\), for \(t\in M\), is a \emph{Hopf surface} generalizing our example (\ref{eqn:hopf-1}) above.  Again we have that \(X_t\) is diffeomorphic to \(S^1\times S^3\) for all \(t\in M\).

Kodaira and Spencer proved the following important theorem in \cite{kodaira-spencer1958}.
\bthm
Let \(t, t'\in M\), then \(X_t\) is biholomorphic to \(X_{t'}\) if and only if there is a \(u\in GL(2,\BC)\) such that 
\[
t'=utu^{-1}.
\]
\ethm
In other words, the Hopf surfaces are biholomorphic if and only if the transformations \(t\) and \(t'\) correspond to a linear change of variables in \(\BC^2\).

Using this theorem we see that there are naturally three classes of distinct Hopf surfaces. Namely, let
\bean
M_0&=&\{t\in M: t =\left(\ba{cc}
								\a&0\\
								0&\d\ea\right), \a\ne \d\},\\
M_1&=& \{t\in M: t=\left(\ba{cc}\a&0\\
									0&\a\ea\right)\},\\
M_2&=&\{t\in M: t=\left(\ba{cc}\a&1\\
								0&\a\ea\right)\}.
\eean

We note that these are all biholomorphically distinct Hopf surfaces, with the single exception that if \(t\in M_0\),
then changing the \emph{order} of the two distinct eigenvalues in the matrix \(t\) yields a biholomorphically equivalent Hopf surface. This follows from the identity
\[
\left(\ba{cc}\a&0\\
				0&\d\ea\right) = \left(\ba{cc}0&1\\
				1&0\ea\right)\left(\ba{cc}\d&0\\
				0&\a\ea\right)\left(\ba{cc}0&1\\
				1&0\ea\right)^{-1}
\]
In a sufficiently small neighborhood \(U\) of any point \(t_0\in M_0\), all the Hopf surfaces \(X_t\) will be biholomorphically distinct.  In fact Kodaira and Spencer show that \(M_0\) is a two-dimensional moduli space for this class of Hopf surfaces.  Namely, the Kodaira-Spencer mapping
\[
\r_t:T_t(M_0)\rightarrow H^1(X_t,\Th_t),
\]
is an isomorphism, where here \(\Th_t\) is the sheaf of holomorphic vector fields on \(X_t\).  In particular, 
\[
\dim T_t(M_0) = \dim H^1(X_t,\Th_1)=2.
\]
The deformation theory for \(M_1\) and \(M_2\) is somewhat more complex (see \cite{kodaira-spencer1958}), but this is not necessary for our purposes here, as we will see below.

We now have our  basic theorems concerning the transcendence degrees of Hopf surfaces. Let \(m\) and \(n\) be nonzero integers and define
\[
Z_{m.n}:= \{t=\left(\ba{cc}\a&0\\
									0&\d\ea\right):\a^m=\d^n\}.
\]
We see that, for each \(m\) and \(n\), \(Z_{m,n}\) is a nonsingular hypersurface in \(M_0\) of complex dimension one.
\bthm
\label{thm:integrality}
If 
\[
t \in Z_{m,n}, M_1,\;\textrm{or}\; M_2\]
then
\[
\deg X_t=1.
\]
\ethm
\pf
First we note that for \(t\in M\), \(\deg X_t\le 1\).  Namely, if \(\deg X_t=2\), it follows from a theorem of Kodaira \cite{kodaira1964} that \(X_t\) would be projective algebraic and hence Kähler.  But this would imply that \(b_1(X_t)\) would be even, but that is not the case since \(X_t\) is diffeomorphic to \(S^1\times S^3\). 

To show that \(\deg X_t=1\) in each of these three cases, it suffices to find a nonconstant meromorphic function \(f\) on \(\BC^2\) which satisfies
\[
f(tz)=f(z).
\]

Suppose that \(t\in Z_{m,n}\), then let 
\[
f(z)=f(z_1,z_2) = \frac{z_1^m}{z_2^n}
\]
then
\bean
f(tz)&=&f(\a z_1,\d z_2)\\
&=&\frac{(\a z_1)^m}{(\d z_2)^n} = \frac{\a^m}{\d^n}\frac{z_1^m}{z_2^n}\\
&=&f(z_1,z_2)=f(z).
\eean

If \(t\in M_1\), simply choose 
\[
f(z)=f(z_1,z_2) = \frac{z_1}{z_2},
\]
and argue in the same way.

If \(t\in M_2\), then it is somewhat more complicated but quite straightforward. Choose
\be
\label{eqn:example-2}
f(z_1,z_2) =\frac{z_1+2z_2+c_1}{z_1+2z_2+c_2},
\ee
where
\bean
c_1&=&-\frac{2\a^2+9\a}{3},\\
c_2&=&\frac{2\a^2+15\a}{3(\a-1)}.
\eean
It is easy to verify that
\bean
f(tz)&=&f(\a z_1+z_2,\a z_2)\\
&=&f(z_1,z_2),
\eean
and thus defines a nonconstant meromorphic function on \(X_t\).  We note that in this example (\ref{eqn:example-2}) the meromorphic function \(f\) on \(\BC^2\) depends on the parameter \(\a\), which was not the case in the first two examples here.
q.e.d.

We now have a classical result concerning Hopf surfaces \(X_t\), for \(t\in M_0\) (see, for instance, Shafarevich \cite{shafarevich1997} or Barth et al \cite{barth2004}).
\bthm
\label{thm:barth}
Let \(t\in M_0\), then
\[\deg X_t= 1
\]
if and only if there exist nonzero integers \(m\) and \(n\) such that
\[
\a^m=\d^n.
\]
\ethm
\pf
By Theorem \ref{thm:integrality} we have seen that if \(\a^m=\b^n\), for some integers \(m\) and \(n\), then \(\deg X_t=1\). We need to show the converse. Suppose that 
\[
\a^m\ne \b^n,
\]
for any integers \(m\) and \(n\).  If \(\deg X_t=1\), then there is a meromorphic function \(f\) on \(W\) such that
\[
f(tz)=f(z).
\]
By Hartogs' theorem for meromorphic functions due to E. E. Levi \cite{levi1910}, \(f\) extends as a meromorphic function to \(\BC^2\) (see, e.g., the classical monographs by Osgood \cite{Osgood1929} and Behnke and Thullen \cite{behnke-thullen1934} for a discussion of this theorem).  We can assume that \(f\) has either a zero or a pole at the origin.  Namely, if \(f\) has a pole st \(z=0\), then this \(f\) suffices.  If \(f\) does not have a pole, then simply replace \(f\) by \(f(z)-f(0)\), and this will also satifsfy our requirement.

By an even older theorem of Poincaré from 1883 \cite{poincare1883}, the function \(f\) can be expressed as the quotient of two holomorphic functions \(g(z)\) and \(h(z)\) on \(\BC^2\),
\[
f(z)=\frac{g(z)}{h(z}.
\]

Let 
\bean
g(z)&=&\sum_{ij} a_{ij}z_1^i z_2^j,\\
h(z)&=&\sum_{kl} b_{kl} z_1 ^k z_2^l,
\eean
be power series representations for \(g\) and \(h\) on \(\BC^2\), then we must have
\[
\frac{g(tz)}{h(tz)}=\frac{g(z)}{h(z)},
\]
which gives
\be
\label{eqn:product}
\left(\sum_{ij}a_{ij}\a^i\d^j z_1^i z_2^j\right)\left(\sum_{kl} b_{kl} z_1^k z_2^l\right) =
\left(\sum_{ij}a_{ij} z_1^i z_2^j\right)\left(\sum_{kl} b_{kl} \a^k\d^lz_1^k z_2^l\right).
\ee

Let \(\mu\) be the order of \(g\) and let \(\nu\) be the order of \(h\), and we have that \(\mu+\nu\ge 1\). Consider the homogenous terms in (\ref{eqn:product}) of lowest order and we have
\be
\label{eqn:product2}
\left(\sum_{i+j=\mu}a_{ij}\a^i\d^j z_1^i z_2^j\right)\left(\sum_{k+l=\nu} b_{kl} z_1^k z_2^l\right) =
\left(\sum_{i+k=\mu}a_{ij} z_1^i z_2^j\right)\left(\sum_{k+l=\nu} b_{kl} \a^k\d^lz_1^k z_2^l\right).
\ee

In this equation, at least one \(a_{ij}\ne 0\) and at least one \(b_{kl}\ne 0\). Let \(a_{{i_0}{j_0}}\) be the unique nonzero coefficient with the property that 
\[
a_{ij}= 0, \;\textrm{for}\; 0\le i<i_0.
\]

Similarly, let \(b_{{k_0}{l_0}}\) be the unique nonzero coefficient such that
\[
b_{kl} = 0,\;\textrm{for}\; 0\le k<k_0.
\]
Then we can rewrite (\ref{eqn:product2}) as
\be
\label{eqn:product3}
\ba{l}
(a_{\mu0}\a^\mu z_1^\mu +\cdots+a_{{i_0}{j_0}}\a^{{i_0}}\d^{j_0}z_1^{i_0}z_2^{j_0})(b_{\nu0} z_1^\nu +\cdots+b_{{k_0}{l_0}}z_1^{k_0}z_2^{l_0})\\
\hspace{.5in}= (a_{\mu0}z_1^\mu +\cdots+a_{{i_0}{j_0}}z_1^{i_0}z_2^{j_0})(b_{\nu0}\a^\nu z_1^\nu +\cdots+b_{{k_0}{l_0}}\a^{{k_0}}\d^{l_0}z_1^{k_0}z_2^{l_0}).
\ea
\ee
We see by expanding these products that the two terms on the left side of the equation
\[
a_{\mu 0}b_{\nu 0}\a^\mu z_1^{\mu+\nu}, \;a_{{i_0}{j_0}}b_{{k_0}{l_0}}\a^{i_0}\d^{j_0}z_1^{i_0+k_0}z_2^{j+0+l_0},
\]
are both unique terms of their respective bidegrees, and this is similar for the terms on the right hand side of the equation%
\footnote{I would like to thank Georges Dloussky for pointing out and correcting an error in this part of my argument that appeared in an earlier draft of this paper.}%
.

It follows that we have the equality
\[
a_{i_0 j_0}b_{k_0 l_0}\a^{i_0}\d^{k_0} z_1^{i_0+k_0}z_2^{j_0+l_0} = a_{i_0 j_0}b_{k_0 l_0}\a^{k_0}\d^{l_0} z_1^{i_0+k_0}z_2^{j_0+l_0},
\]
and  both sides of  this equality are nonzero for  \(z_1\) and \(z_2\) both being nonzero.  This implies that
\[
\a^{i_0-k_0}=\d^{j_0-l_0},
\]
and hence we have a contradiction. q.e.d.

As a simple example, take
\[
t=\left(\ba{cc} 3&0\\0&5\ea\right),
\]
and we see by the fundamental theorem of arithmetic that \(3^m\ne 5^n\), for any integers \(m\) and \(n\), and hence 
\[
\deg X_t =0,
\]
in this case.

The following theorem is now  an easy consequence of this result. Let 
\[
Z:=\bigcup_{m,n} Z_{m,n}.
\]
Here, as before, \(m\) and \(n\) range over the nonzero integers.  We then see that \(Z\) is a countable union of hypersurfaces of dimension one in \(M_0\).  
\bthm
The set 
\[
M_0-Z 
\]
is a set of second category in the metric space \(M_0\), and, in particular, is dense in \(M_0\). Moreover, for each \(t\in M_0-Z\), 
\[
\deg X_t=0.
\]
\ethm
Thus we see that  Hopf surfaces also have the property that generically Hopf surfaces have no nonconstant meromorphic functions.

\section{K3 surfaces}
\label{sec:K3}
In the previous two sections we studied two-dimensional complex tori and Hopf surfaces.  Both of these classes of surfaces are defined as quotients of a fixed Euclidean space (\(\BC^2\) or \(\BC^2-\{0\}\)) by discrete transformation groups acting on these spaces.  In contrast, K3 surfaces are defined by abstract cohomological conditions, and they were first formulated by Weil in 1958 \cite{weil1958} in honor of Kummer, Kähler and Kodaira (as well as a tribute to the Himalayan Peak K2).  They have been studied quite intensely since, and a very readable and informative recent survey is in Buchdahl's paper \cite{buchdahl2003}. In this paper we want to show that a dense set of  K3 surfaces, in a suitably defined moduli space, have transcendence degree zero, i.e., have no nonconstant meromorphic functions.

We will use the standard sheaf-theoretical cohomology theory of several complex variables and complex manifold theory (see, e.g., Wells \cite{wells2008} or Griffiths and Harris \cite{griffiths-harris1978}). In particular, we will use the standard invariants for compact complex manifolds.  Let \(X\) be an \(n\)-dimensional  compact complex manifold, then
\bean
b_i&=&\dim H^i(X,\BC), \;\textrm{Betti numbers},\\
b^+&=& \; \textrm{no. of positive eigenvalues of fundamental quadratic form on}\\
&&\;\;H^{n}(X,\BC),\\
h^{p,q}&=& \dim H^{p,q}(X,\O^p), \;\textrm{Hodge numbers}.
\eean
Moreover, we let \(c_i(E)\) denote the Chern classes of a holomorphic vector bundle on \(X\).

There are various equivalent definitions of K3 surfaces, and we choose the following one. Let \(S\) be a surface, and let 
\[
q=h^{0,1} =\dim H^1(S,\SO),
\]
be the \emph{irregularity} of \(S\), and let
\[
p_g=h^{2,0}=\dim H^0(S,\O^2),
\]
be the \emph{geometric genus} of \(S\).
We define a \emph{K3 surface} to be a surface which has a trivial canonical bundle \(K_S=\wedge^2 T^*(S)\), where \(T(S)\) is the tangent bundle to \(S\), and where the irregularity \(q=0\). It follows that the geometric genus for a K3 surface is simply \(p_g=1\).  

Moreover, for any surface \(S\) with \(q=0\), we note that for any holomorphic line bundle \(F\) on \(S\), \(c_1(F) =0\) if and only if \(F\) is trivial. The Chern class of a surface is defined as the Chern class of its tangent bundle, and thus we have
\bean
c_1(S)&=&c_1(T(S)),\\
&=&c_1(\wedge^2T(S)),\\
&=&-c_1(K_S).
\eean
It follows that \(c_1(S)=0\) if and only if the canonical bundle is trivial, which is an alternative definition of a K3 surface. 

Kodaria also shows in \cite{kodaira1964}that for any surface,
\[
q=\left\{\ba{ll}
b_1/2,\; \textrm{if}\; b_1\;\textrm{is even},\\
(b_1-1)/2,\;\textrm{if}\;b_1\;\textrm{is odd}.
\ea
\right.
\]
Thus, the two defining numerical invariants for K3 surfaces are topological invariants.

In Kodaira's paper \cite{kodaira1964} there are various relationships derived between the various cohomological invariants for the case of surfaces, and, in particular for K3 surfaces, using the above definition.  The resulting Betti and Hodge numbers for K3 surfaces are as follows:
\[
\ba{l}
b_1=0,\\
b_2=22,\\
b^+=3,\\
h^{2,0}=h^{0,2}=1,\\
h^{1,1}=20.
\ea
\]

There are many examples of K3 surfaces which are projective algebraic, in particular, any quartic in \(\BP_3\), such as%
\footnote{We denote projective space over the complex numbers \(\BP_m(\BC)\) simply as \(\BP_m\) in this paper.}%
\[
z_0^4+z_1^4+z_2^4 +z_3^4=0.
\]
is a K3 surface (see, e.g., \cite{griffiths-harris1978}). 

As an illustration, we give the simple proof here.  Let 
\[
Q=\{z\in \BP_3: f(z)=0\},
\]
where \(f\) is a homogeneous polynomial of degree four. by the Leftschetz hyperplane section theorem (e.g. \cite{griffiths-harris1978}),
\[
\ba{lllll}
H^0(Q,\BC)&\cong&H^0(\BP_3,\BC)&\cong&\BC\\
H^1(Q,\BC)&\cong&H^1(\BP_3,\BC)&=&0,
\ea
\]
and thus \(b_1(Q)=0\). Since \(Q\) is a projective algebraic manifold,
\[
b_1=h^{1,0}+h^{0,1},
\]
and hence
\[
q(Q)=h^{0,1}=0.
\]
On the other hand, if \([Q]\) is the line bundle on \(\BP_3\) of the divisor defined by the hypersurface \(Q\), then
\[
[Q]=H^4,
\]
where \(H\) is the hyperplane section bundle on \(\BP_3\). The adjunction formula for a hypersurface in a complex manifold (e.g., \cite{griffiths-harris1978}) gives us the canonical bundle of the hypersurface in terms of the divisor defined by the hypersurface and the canonical bundle of the complex manifold. Namely, in this case here we have
\bean
K_Q&=&K_{\BP_3}\otimes [Q]\\
&=&H^{-4}\otimes H^4,\\
&=&H^0,
\eean
which is the trivial bundle. Thus \(Q\) is a K3 surface.

Moreover, Kodaira discusses in great detail K3 surfaces of transcendence degree one which are all elliptic surfaces, which are fibre spaces over \(\BP_1\) with generic fiber being an elliptic curve and with specific classes of singular fibers.

Let now \(S\) be an arbitrary K3 surface, and, as it is an orientable four dimensional compact differentiable manifold, it has a natural quadratic form acting on \(H^2(S,\BZ)\) which can be defined by
\[
A(\xi,\eta)=\int_S \xi\wedge \eta,
\]
where \(\xi\) and \(\eta\) are closed two-forms representing the integral cohomology classes in \(H^2(S,\BZ)\), noting that torsion plays no role in this quadratic form. The quadratic form \(A\) can be represented by a skew-symmetric integer-valued matrix with signature (3,19), i.e., 3 positive eigenvalues and 19 negative eigenvalues.

Let \(\g_j, j=1,\ldots,22,\) be two-cycles in \(S\) which generate \(H_2(S,\BZ)\), and let \(\eta_i\) be closed two-forms on \(S\) which are dual to \(\{\g_j\}\), i.e.,
\[
\int_{\g_j} \eta_i=\d_{ij}.
\]
This is all independent of the complex structure on \(S\).

Since the canonical bundle \(K_S\) is trivial, then there is a nowhere vanishing holomorphic two form \(\psi\) on \(S\). Moreover, if \(\tilde \psi\) is any other nowhere vanishing two-form on \(S\), then
\[
\tilde \psi =c\psi,
\]
where \(c\) is a constant since,
\[
\dim H^0(S,\O^2)=1.
\]
We now define
\[
\l(S,\psi) =(\l_1,\l_2,\ldots,\l_{22})\in \BC^{22},
\]
by
\[
\l_j = \int_{\g_j} \psi.
\]
Since  any other nowhere vanishing holomorphic two-form differs from \(\psi\) by a constant, we see that \(\l(S,\psi)\) defines a mapping
\[
S\stackrel{\l}{\mapsto} \BP_{21},
\]
which maps \(S\) to a point in the projective space \(\BP_{21}\) which is independent of which two-form \(\psi\) is used.  This is called the \emph{period mapping} for K3 surfaces.

As Kodaira points out in \cite{kodaira1964}, this period mapping which assigns to any K3 surface \(S\) with the given underlying differentiable manifold structure to a point in \(\BP_{21}\) and some of its properties were introduced by Andreotti and Weil. It will be an essential tool in in our investigation of K3 surfaces which are of degree zero, as we shall see in the following paragraphs.

Let now \(\psi\) be a holomorphic two-form on \(S\). It is also a closed two-form. Namely,
\[
d\psi=\partial \psi +\overline{\partial} \psi,
\]
and \(\partial \psi\) is zero by type and \(\overline{\partial}\psi=0\), since \(\psi\) is holomorphic. Thus there are coefficients \((\xi_1,\ldots,\xi_{22})\) such that
\[
\psi=\sum_{i=1}^{22} \xi_i \eta_i + d\f,
\]
where \(\f\) is a one-form on \(S\).  It follows that
\bean
\int_{\g_j}\psi&=&\int_{\g_j} \sum_i \xi_i \eta_i +0,\\
&=&\sum_i \xi_i\int_{\g_j} \eta_i =\xi_j,
\eean
and thus \(\xi_i=\l_i\).  Moreover,
\bean
0&=&\int \psi \wedge \psi,\\
&=&\int_S \left(\sum_i \l_i \eta_i\right) \wedge \left(\sum_j \l_j \eta_j\right),\\
&=&\sum_{i,j}\l_i\l_j\int_S\eta_i\wedge\eta_j,\\
&=&\sum_{i,j} a_{ij}\l_i \l_j.
\eean
Thus \(\l(S)\) is a point on the quadric 
\[
Q:=\{[z_1,\ldots,z_{22}]\in P_{21}: \sum_{i,j} a_{ij}z_i z_j = 0\}.
\]

An important question is: how many K3 surfaces are there?  This has been investigated substantively over the past decades as is outlined in Buchdahl's paper \cite{buchdahl2003}, but for our purposes we will use the local deformation setting formulated and used by Kodaira in his 1964 paper \cite{kodaira1964} using the well known deformation theory of Kodaira and Spencer \cite{kodaira-spencer1958}. We formulate this as the following theorem, and we will outline Kodaira's proof from \cite{kodaira1964}, as it is relatively short and quite informative.  Let \(\Th_S\) be the sheaf of holomorphic vector fields on \(S\). 
\bthm
\label{thm:KNS}
Let \(S\) be a K3 surface, then there is a complex-analytic family \((\SF,\pi,M)\), where \(\SF\stackrel{\pi}{\rightarrow}M\), and \(M\) is an open neighborhood of \(0\in\BC^{20}\), and where
\[
\SF_t:=\pi^{-1}(t), t\in M,
\]
are K3 surfaces diffeomorphic to \(S\), and
\[
S=\SF_{0}=\pi^{-1}(0).
\]
Moreover, the Kodaira-Spencer mapping
\[
\r_0:T(M)_0\rightarrow H^1(S,\Th_S),
\]
is an isomorphism and the family is complete and effectively parametrized.
\ethm
\pf
Let \(\psi\) be a nowhere vanishing holomorphic two-form on \(S\), then it is easy to see that
\be
\label{eqn:isom}
\Th_S\cong \O_S^2.
\ee
Namely, in any local coordinate system \((z_1,z_2)\) on \(S\), we have
\[
\psi =\frac{1}{2}\sum_{\a,\b}\psi_{a\b}dz_\a\wedge dz_{\b},\; \psi_{\a\b}=-\psi_{\b\a}.
\]
If 
\[
v=\sum_\a v_\a\pd{}{z_{\a}},
\]
is a locally defined holomorphic vector field on \(S\), then the mapping
\[
\sum_\a v_\a \pd{}{z_\a}\mapsto \sum_\b v_\a\psi_{\a\b} dz_\b,
\]
induces the required isomorphism (\ref{eqn:isom}).

Kodaira, Nirenberg, and Spencer \cite{kodaira-nirenberg-spencer1958} showed that if a compact complex manifold \(X\) satisfies \(H^2(X,\Th_X)=0\), then there is a complex-analytic family \((\SF,\pi,M)\) with 
\[
\r_0:T_0(M)\rightarrow H^1(X,\Th_X),
\]
being an isomorphism and \(X=\pi^{-1}(0)\), and thus
\[
\dim M=\dim H^1(X,\Th_X).
\]

So in our case, by duality,
\[
h^{2,1}=h^{0,1} = 0,
\]
and hence,
\[
H^2(S,\Th_S)= H^2(S,\O^1)=0,
\]
so the hypothesis of the Kodaira-Nirenberg-Spencer theorem is satisfied., 
and 
\[
\dim H^1(S,\Th_S)=\dim H^1(S,\O^1) =20.
\]
q.e.d.

We note that in the family \(\SF \stackrel{\pi}{\rightarrow} M\) with \(\SF_0=S\), all of the deformations of \(S\) are also K3 surfaces. This follows since \(q\) and \(c_1(S)\) are topological invariants, as noted earlier.

Since each \(\SF_t\) is a K3 surface, then it follows that 
\[
\dim H^1(\SF_t,\Th_t)= 20.
\]
It follows from the stability theorem of Kodaira and Spencer \cite{kodaira-spencer1960} and a completeness theorem of Kodaira and Spencer \cite{kodaira-spencer1958a} that the family is effectively parametrized and complete.
q.e.d.

We now have the following theorem of Kodaira \cite{kodaira1964} which he attributes to Andreotti and Weil. This is often known as a local Torelli theorem in this context of K3 surfaces. We will omit any summary of the proof here.  Essentially, this theorem represents the local moduli parameter space \(M\) in Theorem \ref{thm:KNS} in a quite specific form as an open subset of the quadric \(Q\subset \BP_{21} \) . Here, as before, \(S\) is a given K3 surface and \((\SF,\pi, M)\) is the complex analytic family given in Theorem \ref{thm:KNS}.
\bthm
\label{thm:andreotti-weil}
Let 
\[
\L:M\rightarrow Q\subset \BP_{21}
\]
be defined by
\[
\L(t)=\l(\SF_t),
\]
and let \[
p_0:=\L(0) = \l(S) \in Q.
\]
Then there is an open spherical neighborhood of \(0\in M\) and a neighborhood \(W\subset Q\) of \(p_0\) such that \(\L\) maps \(U\) biholomorphically onto \(W\).
\ethm

We remark that the fibres of the family \(\SF\stackrel{\pi}{\rightarrow} M\) as in Theorem \ref{thm:KNS} are all diffeomorphic, and for \(t\) and \(t'\) in a sufficiently small neighborhood of \(0\in M\), \(\SF_t\) and \(\SF_{t'}\) are biholomorphic if and only if \(t=t'\), i.e., the fibers are locally biholomorphically distinct.

 This theorem allows us to use the coordinate geometry of \(Q \subset \BP_{21}\) to investigate the deformations of \(S\).
 
 Kodaira uses this geometric setting to show that there exist projective algebraic K3 surfaces arbitrarily close to any given K3 surface \(S\) (in the sense of the deformation theory above). We want to now show that arbitrarily close to a given K3 surface \(S\), there are K3 surfaces of  degree zero.
 
 We first have the following Lemma%
\footnote{Thanks to Nicholas Buchdahl for suggesting this lemma and its proof to me.}%
which gives a simple criterion for the existence of nonconstant meromorphic functions.
\blem
Let \(X\) be a compact complex manifold, then if \(f\) is a nonconstant meromorphic function on \(X\), then there is a nontrivial holomorphic line bundle on \(X\).
\elem
\pf Suppose \(f\) is a nonconstant meromorphic function on \(X\), then if \((f)\) is the divisor associated to the meromorphic function \(f\), then \((f)\) is the difference of two effective divisors,
\[
(f)=D_+ - D_-,
\]
where \(D+\) corresponds to the zero set of \(f\) and \(D_-\) corresponds to the polar set of \(f\), both of which are subvarieties of \(X\).  Thus, either of the two divisors \(D_+\) or \(D_-\) must then give rise to a nontrivial holomorphic line bundle on \(X\). 	q.e.d.

Now we have a lemma due to Kodaira \cite{kodaira1964}, which helps characterize nontrivial holomorphic line bundles on a K3 surface in terms of the bilinear form \(A\) on \(\BC^{22}\).  The proof is not difficult, and we refer the reader to Kodaira's paper. Again \(S\) is a fixed K3 surface.
\blem
A cohomology class \(c\in H^2(S,\BZ)\) is the Chern class of a holomorphic line bundle \(F\) over \(S\), if and only if the point
\[
m=(m_1,\ldots,m_{22}), \;m_j=\int_{\g_j}  c,
\]
satisfies the linear equation
\[
A(\l,m)=0,
\]
where \(\l= \l(S)\).
\elem
We note that the line bundle \(F\) is nontrivial if and only if \(m\ne0\), since the irregularity \(q(S)=0\).

We can now formulate our fundamental result concerning K3 surfaces with no nonconstant meromorphic functions. We let \(S\) be an arbitrary K3 surface, and let \(\SF\stackrel{\pi}{\rightarrow} M\) be the local deformation space with \(\SF_0=S\), where \(t\) is the local parameter in \(M\), as given in Theorem \ref{thm:KNS}.
\bthm
For a sufficiently small neighborhood  \(U\) of \(0\) in \(M\), there is a set \(\O\subset U\) of second category, which, in particular, is dense in \(U\), such that if \(t\in \O\), then
\[
\deg \SF_t = 0.\]
\ethm
\pf
We use Theorem \ref{thm:andreotti-weil} to biholomorphically map a sufficiently small neighborhood \(U\) of \(0\) onto a neighborhood \(W\subset Q\) of \(p_0 = \L(0)\).  For any \(\L (t)\), for \(t\in U\), \(\SF_t\) has a nontrivial holomorphic line bundle \(
F\) if and and only if
\[
\int_{\g_j} c(F) ,
\]
is a nonzero integer for some \(j, j=1,\ldots,22\). Let now
\[
m=(m_1,\ldots,m_{22})\in \BZ^{22}-\{0\},
\]
and we let
\[
\tilde Z_m:= \{z\in \BC^{22}:A(z,m)=0\}.
\]
Each such \(\tilde Z_m\) defines a linear hyperplane in \(\BP_{21}\), and we let 
\[
Z_m:=\tilde Z_m\cap W, 
\]
which is a complex subvariety of \(W\) with \(\dim Z_m\le 19\), and, of course, the intersection could well be an empty set for some \(m\).  In any event, \(Z_m\) is a closed subset of \(W\) with no interior points.  

We let
\[
Z:= \bigcup_{m\in\BZ^{22}-\{0\}} Z_m,
\]
and this is a countable union of closed subsets of \(W\) with no interior points, and hence
\[
\O:=W-Z,
\]
is a set of second category in \(W\) which is dense in \(W\). As we saw earlier, any point \(p\in W\) with \(\deg \SF_{\L^{-1}(p)} \ge 1\) must be a point in \(Z_m\) for some \(m\ne0\).  It follows that if \(p\in \O\),
\[
\deg \SF_{\L^{-1}(p)} =0.
\]
q.e.d

\bibliography{references}

\begin{thebibliography}{10}

\bibitem{barth2004}
Wolf~P. Barth, Claus~Klaus Hulek, Christ A.~M. Peters, and Antonius van~de Ven.
\newblock {\em Compact Complex Surfaces}.
\newblock Springer, Berlin, 2004.

\bibitem{behnke-thullen1934}
H.~Behnke and P.~Thullen.
\newblock {\em {Theorie der Funktionen mehrerer komplexer Veränderlichen}}.
\newblock Springer-Verlag, Heidelberg, 1934.

\bibitem{buchdahl2003}
Nicholas Buchdahl.
\newblock {Compact Kähler surfaces with trivial canonical bundle }.
\newblock {\em Annals of Global Analysis and Geometry}, 23:189--204, 2003.

\bibitem{conforto1956}
Fabio Conforto.
\newblock {\em Abelsche Funktionen und Algebraische Geometrie}.
\newblock Springer-Verlag, Berlin, 1956.

\bibitem{grauert1962}
Hans Grauert.
\newblock {Über Modifikationen und exzeptionelle analytische Mengen}.
\newblock {\em Math. Annalen}, 146:331--368, 1962.

\bibitem{griffiths-harris1978}
Phillip Griffiths and Joseph Harris.
\newblock {\em Principles of Algebraic Geometry}.
\newblock John Wiley \& Sons, New York, 1978.

\bibitem{hartogs1906}
Fritz Hartogs.
\newblock { Einige Folgerungen aus der Cauchyschen Integralformel bei
  Funktionen mehrerer Veränderlichen}.
\newblock {\em Sitzungsberichte der Königlich Bayerischen Akademie der
  Wissenschaften zu München, Mathematisch-Physikalische Klasse}, 36:223--242,
  1906.

\bibitem{hopf1948}
Heinz Hopf.
\newblock {Zur Topologie der komplexen Mannigfaltigkeiten}.
\newblock In {\em Studies and Essays Presented to R. Courant on his 60th
  Birthday}, pages 167--185. Interscience Publishers, Inc., New York, 1948.

\bibitem{kodaira1964}
K.~Kodaira.
\newblock {On the structure of compact complex surfaces. I}.
\newblock {\em Amer. J. Math.}, 86:751--798, 1964.

\bibitem{kodaira-nirenberg-spencer1958}
K.~Kodaira, L.~Nirenberg, and D.~C. Spencer.
\newblock On the existence of complex analytic structures.
\newblock {\em Annals of Math.}, 68:450--459, 1958.

\bibitem{kodaira-spencer1958}
K.~Kodaira and D.~C. Spencer.
\newblock {On deformations of complex-analytic structures I, II}.
\newblock {\em Annals Math.}, 67:328--466, 1958.
\newblock Published as two consecutive papers in the same volume of the Annals
  of Mathematics.

\bibitem{kodaira-spencer1958a}
K.~Kodaira and D.~S. Spencer.
\newblock A theorem of completeness for complex analytic fibre spaces.
\newblock {\em Acta Math.}, 100:281--294, 1958.

\bibitem{kodaira-spencer1960}
K.~Kodaira and D.~S. Spencer.
\newblock {On deformations of complex analytic structures, III. Stability
  theorems for complex structures}.
\newblock {\em Annals of Math.}, 71:43--76, 1960.

\bibitem{levi1910}
E.~E. Levi.
\newblock Studii sui punti singolari essentiali delle funzioni analitiche di
  due o pi\`{u} variabili.
\newblock {\em Ann. di Mat.}, 17:61--88, 1910.

\bibitem{Osgood1929}
W.~F. Osgood.
\newblock {\em Lehrbuch der Funktionentheorie, Vol 2, Pt. 1}.
\newblock B. G. Teubner, Leipzig, 1929.

\bibitem{poincare1883}
Henri Poincaré.
\newblock Sur les fonctions de deux variables.
\newblock {\em Acta Math.}, 2:97--113, 1883.

\bibitem{riemann1857}
Bernhard Riemann.
\newblock {Theorie der Abel'schen Functionen}.
\newblock In H.~Weber, editor, {\em Gesammelte Mathematische Werke}, pages
  81--135. Teubner, Leipzig, 1876.
\newblock This paper in Riemann's Werke combines four paper from the Journal
  für reine und angewandte Mathematik, Vol. 54, 1857.

\bibitem{shafarevich1997}
Igor~R. Shavarevich.
\newblock {\em Basic Algebraic Geometry 2}.
\newblock Springer-Verlag, Berlin, 1997.

\bibitem{siegel1949}
Carl~L. Siegel.
\newblock {\em Analytic Functions of Several Complex Variables}.
\newblock Institute for Advanced Study, Princeton, NJ, 1949.
\newblock Lecture Notes 1948--1949, Notes by P. T. Bateman.

\bibitem{siegel1955}
Carl~L. Siegel.
\newblock {Meromorphe Funktionen auf kompakten analytischen
  Mannigfaltigkeiten}.
\newblock {\em Nachrichten der Akademie der Wissenschaften un Göttingen.
  Mathematische-physikalische Klasse}, 4:71--77, 1955.
\newblock Also in Carl Ludwig Siegel Gesamelte Abhandlungen, Band III,
  Springer-Verlag, 1966, pp. 216--222.

\bibitem{weil1958}
André Weil.
\newblock Final report on contract af 18 (608)-57, 1958.
\newblock In {\em Oeurves Scientifique, Vol. II}. Springer, New York, 1980.

\bibitem{wells2008}
Raymond~O. Wells, Jr.
\newblock {\em Differential Analysis on Complex Manifolds}.
\newblock Springer-Verlag, New York, 2008.
\newblock 3rd Edition.

\bibitem{weyl1913}
Hermann Weyl.
\newblock {\em Die Idee der Riemannschen Fläche}.
\newblock B. G. Teubner, Leipzig, 1913.

\end{thebibliography}
\bibliographystyle{plain}%

\end{document}